\documentclass[12pt]{article}
\usepackage{latexsym, amscd, amsmath, graphics, rotating, array}
\usepackage{amssymb,amsmath}
\makeindex

\title{Zeta functions of groups and enumeration in Bruhat-Tits buildings}
\author{Christopher Voll\thanks{Mathematical Institute, 24-29 St. Giles', Oxford, OX1 3LB, United Kingdom}
\thanks{Supported by an EPSRC Postdoctoral Research Fellowship in Mathematics.}}

\date{30 September 2002}

\newenvironment{remark}{\bigskip\noindent\textsl{Remark.}\rm}

\newenvironment{acknowledgements}{\bigskip\noindent\textsl{Acknowledgements.}\rm}

\newenvironment{example}{\bigskip\noindent\textbf{Example. }\,\rm}
{\bigskip}

\newenvironment{proof}{\bigskip\noindent\textbf{Proof.}\rm}{\hfill$\Box$ \\ \hspace{0.1cm}}

\newenvironment{proofnodot}{\bigskip\noindent\textbf{Proof}\rm}{\hfill$\Box$ \\ \hspace{0.1cm}}

\newtheorem{lemma}{Lemma}
\newtheorem{theorem}{Theorem}

\newtheorem{corollary}{Corollary}
\newtheorem{proposition}{Proposition}
\newtheorem{definition}{Definition}
\newtheorem{observation}{Observation}

\def \N {\ensuremath{\mathbb{N}}}
\def \coneN {N}

\def \oK{\mathfrak{o}_K}

\def \Q {\mathbb{Q}}
\def \u {\"u}

\def \Woffoff {A_{\mbox{\small off/off }}(p,T)}
\def \Wsmptoff {\ensuremath{A_{\mbox{\small sm.pt./off}}(p,T)}}

\def \T {\ensuremath{\mathfrak{T} }}

\def \Qp {\mathbb{Q}_p}
\def \Z {\mathbb{Z}}
\def \z {\zeta^\triangleleft}
\def \Zp  {\mathbb{Z}_p}

\def \Fp  {\mathbb{F}_p}

\def \SF {\ensuremath{\mathcal{S}_F}}
\def \Fpthree {\ensuremath{\mathcal{F}(p,3)}}

\def \ApTSF {\ensuremath{A(p,T,F)}}

\def \PtwooverFp {\mathbb{P}^2(\mathbb{F}_p)}

\def \SlnQp {\ensuremath{Sl_n({\mathbb{Q}_p})}}

\def \L {\ensuremath{\Lambda}}

\def \Lcen {\ensuremath{\Lambda_{der}}} 
\def \Lab {\ensuremath{\Lambda_{ab}}}
\def \Mcen {M_{der}}
\def \Mab  {M_{ab}}
\def \LieL {\ensuremath{\mathfrak{L}}}

\def \D {\ensuremath{\Delta}}

\def \fFp {\ensuremath{n_{f,p}}}
\def \a {\ensuremath{\alpha}}

\def \F23 {\ensuremath{F_{2,3}}}
\def \Dstar {\ensuremath{\mathfrak{D}^*}}

\def \pthreetwo {\binom{3}{2}_p}
\def \ptwoone {\binom{2}{1}_p}

\newcommand{\verylongpage}{\enlargethispage*{2.3cm}}

\begin{document}
\maketitle
\begin{abstract}
We introduce a new method to calculate local normal zeta functions of
finitely generated, torsion-free nilpotent groups, \T-groups in
short. It is based on an enumeration of vertices in the Bruhat-Tits
building for~\SlnQp. It enables us to give explicit computations for
\T-groups of class $2$ with small centres and to derive local
functional equations. Examples include formulae for non-uniform normal
zeta functions. 
\end{abstract}


\section{Introduction and results}\label{section1}
A finitely generated group has only a finite number of subgroups of each finite index. With this simple observation Grunewald, Segal and Smith opened their seminal paper~\cite{GSS/88} by which they initiated the study of zeta functions of torsion-free finitely generated nilpotent groups, \T-groups in short. To a \T-group $G$ they associate a Dirichlet series
$$\z_{G}(s):=\sum_{n=0}^\infty a_n n^{-s},$$
where $$a_n:=|\{H\triangleleft G|\;|G:H|=n\}|$$ and $s$ is a complex variable. This series is called the {\sl normal zeta function} of~$G$. It decomposes as an Euler-product
$$\z_{G}(s)=\prod_{\text{$p$ prime}} \z_{G,p}(s),$$
of {\sl local} normal zeta functions $$\z_{G,p}(s):=\sum_{n=0}^\infty a_{p^n} p^{-ns}.$$ One of the main results of~\cite{GSS/88} is the {\sl rationality} of the local factors. However, the local zeta function's dependence on the prime $p$ remained mysterious and stayed in the focus of research in the subject. The question was recently linked to the classical problem of counting points on varieties mod $p$ by du Sautoy and Grunewald~(\cite{duSG/00}). In the present paper we introduce a new method to compute local normal zeta functions of nilpotent groups which exploits the combinatorial geometry of the Bruhat-Tits building for $\SlnQp$. It allows us to prove
\begin{theorem} \label{main theorem}
Let $G$ be a \T-group of nilpotency class $2$ with derived group $G'$. Assume that the associated Lie ring $\LieL(G):=G/G'\oplus G'$ has a presentation as follows:
\begin{equation*}
\LieL(G)=\langle x_1,\dots,x_{2r},y_1,y_2,y_3|\;(x_i,x_j)=M({\bf y})_{ij}\rangle,
\end{equation*}
where $M({\bf y})=\left(\begin{array}{cc}
                            0 &R({\bf y})\\
                         -R({\bf y})^t&0 
                        \end{array}\right)$ is a matrix of $\Z$-linear forms in ${\bf y}=(y_1,y_2,y_3)$, $r\geq2$, and all other Lie-brackets are understood to be trivial. Assume that the curve $$C=\{{\bf y}\in \mathbb{P}^2(\mathbb{Q})|\;\det(R({\bf y}))=0\}$$ is smooth. Set $$|C(\Fp)|:=|\{{\bf y}\in \mathbb{P}^2(\Fp)|\;\det(R({\bf y}))=0\}|.$$ Then for almost all primes $p$
$$\z_{G,p}(s)=W_1(p,p^{-s})+|C(\Fp)|W_2(p,p^{-s})$$ for explicitly determined rational functions $W_i(X,Y)\in\Q(X,Y)$, $i=1,2$.

\end{theorem}

By a result of Beauville~(\cite{Beauville/00}, Proposition~3.1) every
smooth plane curve over~$\Q$ can be defined by a linear
determinant. Recall that the {\sl Hirsch length} of a~$\T$-group~$G$ is
the number of infinite cyclic factors in a decomposition series for~$G$. We deduce

\begin{corollary} Given a smooth plane curve $C$ over $\Q$ of degree $r\geq2$, there is a \T-group $G=G_C$ of Hirsch length $2r+3$ and rational functions $W_i(X,Y)\in\Q(X,Y)$, $i=1,2$ such that for almost all primes $p$
\begin{equation*}
\z_{G_C,p}(s)=W_1(p,p^{-s})+|C(\Fp)|W_2(p,p^{-s}). 
\end{equation*}
\end{corollary} 

The explicit expressions for the $W_i(X,Y)$ produced in the proof of Theorem~\ref{main theorem} allow us to read off a functional equation of the type 
\begin{equation}
W_i(X^{-1},Y^{-1})=(-1)^{n_i}X^{a_i}Y^{b_i}W_i(X,Y),\;n_i,a_i,b_i\in\Z.
\label{combinatorial functional equation}
\end{equation}
The {\sl rationality} of the Weil zeta function for an irreducible variety $V$ over a finite field $\Fp$ implies that the function 
\begin{eqnarray*}
\N_{\geq0}&\rightarrow&\N \\
e&\mapsto&|V(\mathbb{F}_{p^e})|:=|\{\mathbb{F}_{p^e}-\mbox{rational points of }V\}|
\end{eqnarray*}
 has a unique extension to $\Z$. In particular, the symbol $|V(\mathbb{F}_{p^{-1}})|$ is well-defined. The {\sl functional equation} satisfied by the Weil zeta function for a non-singular, absolutely irreducible projective variety $V$ over $\Fp$ implies the relation
\begin{equation}
|V(\mathbb{F}_{p^{-e}})|=p^{-en}|V(\mathbb{F}_{p^{e}})|.\label{weil functional equation}
\end{equation}
Combining the functional equations~(\ref{combinatorial functional equation}) and~(\ref{weil functional equation}) will allow us to deduce

\begin{corollary}\label{corollary to main theorem: fun. eq.} Assume further that the curve $C$ in Theorem~\ref{main theorem} is absolutely irreducible over $\Q$. Then for almost all primes $p$
\begin{equation}
\z_{G_C,p}(s)|_{p\rightarrow p^{-1}}=-p^{\binom{2r+3}{2}-(4r+3)s}\z_{G_C,p}(s). \label{fun. eq. main theorem}
\end{equation}
\end{corollary}

Theorem~\ref{main theorem} produces a wealth of concrete examples of
normal zeta functions which are not finitely uniform. We call
$\z_{G}(s)$ {\sl finitely uniform} if the primes fall into finitely
many classes on which the local factors are given by one rational
function, i.e. if there are rational functions
$W_1(X,Y),\dots,W_n(X,Y)$ such that for all primes $p$ there is an
$r=r(p)\in\{1,\dots,n\}$ such that $\z_{G,p}(s)=W_r(p,p^{-s})$. We say
that $\z_{G}(s)$ is {\sl uniform} if $n=1$. Du Sautoy was the first to
construct examples of zeta functions which are not finitely uniform
(\cite{duS-ecI/01},\cite{duS-ecII/01}). He constructed Lie rings presented by matrices $M({\bf y})$ as in Theorem~\ref{main theorem} above whose determinants define certain elliptic curves over $\Q$. His analysis established the {\sl existence} of the rational functions $W_i(X,Y)$ and did not reveal the functional equation~(\ref{fun. eq. main theorem}). His Theorem 1.1 and the assertion of Conjecture 5.6 in~\cite{duS-ecII/01} follow from our Theorem~\ref{main theorem} above.

Our second theorem confirms that these examples were indeed in some sense minimal. Our method allows us to describe explicitly the local normal zeta functions of class-$2$-nilpotent \T-groups with derived groups of Hirsch length~$2$. 

\begin{theorem} \label{general rank 2 centres theorem}
Let $G$ be a \T-group of class $2$ with derived group $G'$ of Hirsch length~$2$. Then there are irreducible polynomials $f_1(t),\dots,f_m(t)\in \Q[t]$ and rational functions $W_I(X,Y)$, $I\subseteq \{1,\dots,m\}$ such that for almost all primes~$p$
$$\z_{G,p}(s)=\sum_{I\subseteq\{1,\dots,m\}}c_{p,I}W_{I}(p,p^{-s}),$$
where 
\begin{equation}
c_{p,I}=|\{x\in \mathbb{P}^1(\mathbb{F}_p)| \,f_i(x)\equiv0\mbox{ mod }  p\mbox{ if and only if } i\in I\}|.\label{definition cpI}
\end{equation}
In particular, $\z_{G}(s)$ is finitely uniform. For almost all primes $\z_{G,p}(s)$ satisfies the functional equation
$$\z_{G,p}(s)|_{p\rightarrow p^{-1}}=(-1)^{d+2}p^{\binom{d+2}{2}-(2d+2)s}\z_{G,p}(s).$$
\end{theorem}

The polynomials $f_i$ in Theorem~\ref{general rank 2 centres theorem} are those occurring in the classification of torsion-free {\sl radicable} nilpotent groups of class~$2$ and of finite Hirsch length with centres of Hirsch length~$2$ due to Grunewald and Segal~(\cite{GSegal/84}) on which the proof of Theorem~\ref{general rank 2 centres theorem} relies. We will recall this classification in Section~\ref{section3} to make this paper self-contained. The proof of Theorem~\ref{general rank 2 centres theorem} will also produce an algorithm for the computation of the $W_I(X,Y)$.

In Section~\ref{section2} we will explain the relationship between the combinatorics of buildings and local zeta functions. We shall prove Theorem~\ref{general rank 2 centres theorem} in Section~\ref{section3}. Its proof uses the classification of torsion-free {\sl radicable} nilpotent groups of class~$2$ and of finite Hirsch length with centres of Hirsch length~$2$ by Grunewald and Segal~(\cite{GSegal/84}). In Section~\ref{section4} we will give the proof of Theorem~\ref{main theorem} and its corollaries as well as the explicit formulae for du Sautoy's example. 

We feel that the functional equations for non-uniform zeta functions may well be the most important result in the present paper. ``To find an explanation of the phenomenon in general [$\dots$] is one of the most intriguing open problems in this area'', say du Sautoy and Segal in~\cite{duSSegal/00}.  

\begin{acknowledgements}
This article comprises parts of the author's Cambridge PhD thesis~\cite{Voll/02} which was supported by the Studienstiftung des deutschen Volkes and the Cambridge European Trust. The author should like to thank Marcus du Sautoy for his continuous support and encouragement as supervisor and Fritz Grunewald for numerous invaluable conversations in D\u sseldorf. Thanks are also due to Burt Totaro who directed us to~\cite{Beauville/00}.
\end{acknowledgements}

\section{Zeta functions of groups and the Bruhat-Tits building $\D_n$ for \SlnQp}\label{section2}

In this section we will explain how counting ideals in \T-groups may be interpreted as enumerating certain vertices in the Bruhat-Tits building $\D_n$ for~\SlnQp. In nilpotency class~$2$, we describe how this approach may be used to facilitate effective computations, at least if the derived group is small.

Let $G$ be a \T-group. Let $\LieL=\LieL(G)$ be the $\Q$-Lie algebra associated to $G$ under the Malcev correspondence. We can define the local ideal zeta functions associated to $\LieL$ similarly to the local normal zeta functions associated to $G$. For a prime $p$ we set
$$\z_{\LieL,p}(s):=\sum_{n=0}^\infty b_{p^n}p^{-ns},$$
where $$b_{p^n}:=|\{\Lambda \triangleleft \LieL|\;|\LieL:\Lambda|=p^n\}|.$$
Theorem 4.1 in~\cite{GSS/88} confirms the following:

\begin{proposition}
For almost all primes $p$, 
\begin{equation}
\z_{G,p}(s)=\z_{\LieL(G),p}(s).\label{normal sbgps vs. ideals}
\end{equation}
\end{proposition}

If $G$ is a \T-group of nilpotency class $2$, we may indeed work with the graded Lie ring $\LieL(G)=G/G'\oplus G'$ and equation~(\ref{normal sbgps vs. ideals}) holds for {\sl all} primes $p$ (this is essentially the Remark in~\cite{GSS/88}, p. 206). Ideals of $p$-power index in~$\LieL$ correspond to ideals of finite index in $\LieL_p:=\LieL\otimes_{\Z}\Zp$. Note that if $G$ has Hirsch length $n$ then $\LieL_p\cong \Zp^n$ as $\Zp$-modules for almost all primes $p$. 

That the Bruhat-Tits building $\D_n$ for $\SlnQp$ is a natural place to represent ideals in $\LieL_p$ follows from the following trivial 
\begin{observation}\label{observation on ideals}
A lattice $\Lambda \leq \LieL$ is an ideal if and only if $p^n \Lambda$ is an ideal for all $n\in\Z$.
\end{observation}

\subsection{The Bruhat-Tits building $\D_n$ for $\SlnQp$}

We will recall some features of $\D_n$ which we will need for our
analysis. We refer the reader to~\cite{Serre/80}, \cite{Garrett/97},
\cite{Brown/89} for more background on buildings.

The building $\D_n$ may be viewed as a simplicial complex whose vertices are {\sl homothety classes} of lattices in $\Qp^n$. Recall that lattices $\L$ and $\L'$ are called homothetic if there is an element $x\in\Qp$ such that $\L=x\L'$. We will use square brackets to denote a lattice's (homothety) class. One defines an {\sl incidence relation} on vertices as follows: Classes $X_1$ and $X_2$ are called incident if there are lattices $\L_1$ and $\L_2$ with $X_i=[\L_i], i=1,2$, such that $$p\L_1<\L_2<\L_1.$$
The building $\D_n$ is just the {\sl flag complex} for this incidence geometry: its simplices consist of sets of pairwise incident vertices. For our purposes we are only interested in lattices in $\Qp^n$ which are contained in the {\sl standard lattice} $\Zp^n$. From now on, {\sl lattice} will mean {\sl lattice contained in $\Zp^n$} and {\sl class} will mean {\sl homothety class of lattices contained in~$\Zp^n$}. It is then clear that every class contains a unique ($\leq$-) maximal lattice. We call~$[\Zp^n]$ the {\sl special vertex} or {\sl root vertex}.

Observation~\ref{observation on ideals} amounts to saying that ideals define a simplicial subcomplex of $\D_n$. Before we show how this subcomplex may be described effectively for a class-$2$-nilpotent \T-group, we have to set up some notation.

Every building is glued together from apartments. The apartments of~$\D_n$ are simplicial subcomplexes, isomorphic to Euclidean $(n-1)$-space tessellated by $(n-1)$-dimensional simplices. There is a $1-1$-correspondence between apartments and decompositions of $\Qp^n$ into a direct sum of lines. Any choice of generators of these lines gives a basis for $\Qp^n$, and vertices contained in a fixed apartment correspond to lattices with $p$-power multiples of this basis' elements as elementary divisor basis. 

\begin{figure}\label{figure-sectorbundle}
\begin{center}
\resizebox{7cm}{!}{\includegraphics{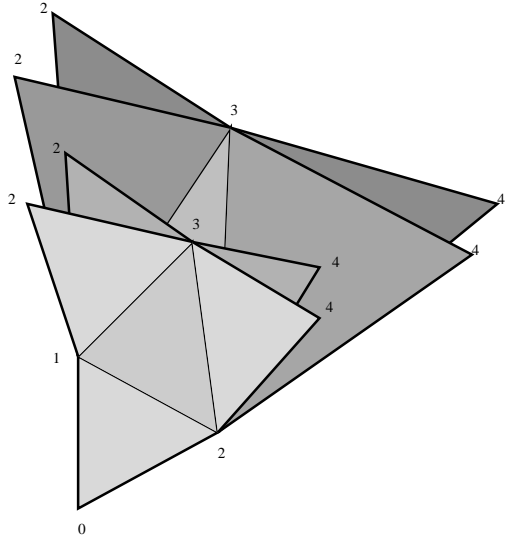}} 
\caption{Initial segment of a sector-family in $\Delta_3$ $(p=2)$. Vertices $[\L]$ are labelled by $w([\L])$ (cf Definition~\ref{definition:new-weight-function} of~$w$).}
\end{center}
\end{figure}

Apartments containing the special vertex will be called {\sl special}. There is a $1-1$-correspondence between special apartments and decompositions of $\Zp^n$ into a direct sum of $\Zp$-modules of rank $1$. 
Each special apartment falls naturally into $|S_n|=n!$ {\sl sectors} with the root as vertex: given an {\sl ordered} $\Zp$-basis $(e_1,\dots,e_n)$ for $\Zp^n$ and a permutation $\sigma\in S_n$, the sector corresponding to $\sigma$ is defined to be the simplicial subcomplex generated by classes of lattices $$\langle p^{r_{1}}e_{\sigma(1)},\dots, p^{r_{n}}e_{\sigma(n)} \rangle,\;r_{1}\geq r_{2} \geq \dots \geq r_{n}\in\N.$$ Note that this decomposition of special apartments into sectors is not disjoint. 

Every $1$-codimensional simplex is contained in $p+1$ maximal simplices and any vertices' {\sl link}\footnote{Recall that the link of a vertex $v$ in a simplicial complex is $\overline{St(v)} - St(v)$, where the vertices' {\sl star} $St(v)$ is the union of all open simplices having $v$ as vertex.} is isomorphic to the {\sl flag complex}~$\mathcal{F}(\Fp^n)=\mathcal{F}(p,n)$, say, whose vertices are the proper subspaces of~$\Fp^n$ together with the obvious incidence relation. 
We say that sectors are {\sl equivalent mod $p$} if they have the first maximal simplex (i.e. the one containing the root) in common. We will call the equivalence classes {\sl sector-families}. They are indexed by complete flags $F$ in $\Fp^d$. A sector's {\sl boundary} is defined to be the simplicial subcomplex generated by the subset of its vertices which are contained in at least two sectors\footnote{ These are clearly the classes of lattices whose elementary divisor type contains at least two repeated $p$-powers.} inequivalent mod $p$. A sector-family's {\sl boundary} is defined to be the union of the boundaries of its sectors. We will use the notation $\mathcal{S}_F$, $\partial\mathcal{S}_F$, $\mathcal{S}^{\circ}_F$ to denote the sector-family indexed by the flag~$F$, its boundary and interior, respectively. Notice that sector-families overlap at their boundaries. 

\begin{remark} Sector-families have a natural interpretation as fibres of maximal simplices under the ``reduction mod $p$''-map $\pi:X^{\infty}\D_n\rightarrow X^0\D_n$ from the (spherical) Tits building at infinity to the (spherical) Tits building ``around the special vertex''. These are just the usual buildings for $Gl_n(\Qp)$ and $Gl_n(\Fp)$, respectively. We should like to thank Linus Kramer for pointing this out to us.
\end{remark}

\subsection{$\D_n$ and \T-groups of class $2$}

From now on, let $\LieL$ be a Lie ring corresponding to a class $2$ nilpotent \T-group, with derived ring $\LieL'$, finitely presented by
\begin{equation}
\LieL=\langle x_1,\dots,x_d,y_1,\dots,y_{d'}|\;(x_i,x_j)=M({\bf y})_{ij}\rangle,\label{Lie ring presentation}
\end{equation}
where $M({\bf y})$ is a matrix of $\Z$-linear forms in $y_1,\dots,y_{d'}$.

Following~\cite{Segal/83}, Chapter~11, Section~C, we define a full alternating bilinear map
\begin{eqnarray}
\phi_{\LieL}:(\LieL/\LieL')\times (\LieL/\LieL')&\rightarrow& \LieL'\label{bilinear map} \\
(a\LieL',b\LieL')&\mapsto&(a,b). \nonumber
\end{eqnarray}
We say that the matrix $M$ represents the map $\phi_{\LieL}$ and refer to the square root of the determinant $M$ as the Pfaffian associated to $\phi_{\LieL}$ or even to $\LieL$. Let~$p$ be a prime number such that $\LieL_p/{\LieL_p'}$ is torsion-free. In this chapter we show how the local normal zeta function $\z_{\LieL,p}(s)$ depends on the bilinear map~(\ref{bilinear map}) and how it can be computed effectively in terms of an integer-valued weight function on the vertices of the Bruhat-Tits building $\D_{d'}$.

A lattice $\L$ in $\LieL_p$ gives rise to a pair $(\Lab,\Lcen)$ of a lattice $$\Lcen:=\L\cap \LieL_p'\leq\LieL_p'$$ in the derived ring and a lattice $$\Lab:=(\L+\LieL_p')/\LieL_p'\leq\LieL_p/\LieL_p'$$ in the abelian quotient. It is easy to see that this map is surjective and $|\LieL_p':\Lcen|^d$ to~$1$. 
The lattice~$\L$ defines an ideal if and only if its associated pair~$(\Lab,\Lcen)$ is {\sl admissible} in the following sense.

\begin{definition}\label{admissible pairs}
A pair $(\Lab,\Lcen)$, $\Lab\subseteq\LieL_p/\LieL_p',\,\Lcen\subseteq\LieL_p'$ is {\rm admissible} if 
\begin{equation}
\forall j \in \{1,\dots,d\}:\;\phi_j\Lambda_{ab} \subseteq \Lcen  \label{latticecondition}
\end{equation}
where $\phi_j$ denotes the linear map $v \mapsto (v,x_j)$ and~$x_j,\,j=1,\dots,d$ are the Lie generators of~$\LieL$. A lattice $\Lab$ is called {\rm admissible for~$\Lcen$} if~(\Lab,\Lcen) form an admissible pair.
\end{definition}

Keeping track of the index of the largest lattice $\Lab$ admissible for a given~$\Lcen$ is sufficient to compute the normal zeta function. More precisely, we put 
$$X(\Lcen)/\Lcen:=Z(\LieL_p/\Lcen)$$ 
and have
\begin{equation}
\z_{\LieL_p}(s)=\z_{\Zp^d}(s)\sum_{\Lcen\subseteq\LieL_p'}|\LieL_p':\Lcen|^{d-s}|\LieL_p:X(\Lcen)|^{-s}. \label{GSS-Lemma6.1}
\end{equation} 

This is essentially Lemma 6.1. in~\cite{GSS/88}. Our aim is to understand the function $\Lcen\mapsto|\LieL_p:X(\Lcen)|$. Let us assume first that $\Lcen$ is maximal in its class and of elementary divisor type $$\nu=(p^{r_1},p^{r_2},\dots,p^{r_{d'-1}},1),\;r_1\geq r_2\geq\dots\geq r_{d'-1}\in\N.$$  If $\Mcen$ is a matrix whose rows contain the coordinates of generators for $\Lcen$ over $\Zp$, written with respect to the chosen $\Zp$-basis for $\LieL_p$, we will write $\Lcen=\Zp^{d'}\cdot \Mcen$, which is standard notation. As $Sl_{d'}(\Zp)$ acts transitively on lattices with a fixed elementary divisor type there exists $\alpha \in Sl_{d'}(\Zp)$ such that 
$$\Zp^{d'}\cdot\Mcen\alpha=\Zp^{d'}\cdot\mbox{diag}(p^{r_1},p^{r_2},\dots,p^{r_{d'-1}},1),$$
Note that $\alpha$ is unique only up to multiplication from the right by an element of~$G_{\nu}$, the stabilizer in~$Sl_{d'}(\Zp)$ of the lattice spanned by the diagonal matrix. We write $\alpha^i$ for the $i$-th column of the matrix~$\alpha$ and~$\alpha_i^j$ for its~$ij$-th entry. Let $C_j$ be the matrices of the linear maps~$\phi_j$ (``Lie-bracketing with the~$j$-th Lie-generator'') for $j=1,\dots,d$. The condition~(\ref{latticecondition}) for a lattice~$\Lab=\Zp^d\cdot\Mab$ to be admissible for the lattice~$\Lcen=\Zp^{d'}\cdot\Mcen$ may be rephrased as
\begin{eqnarray} 
\forall j \in \{1,\dots, d\}:\;\Zp^{d'}\cdot\Mab C_j & \subseteq & \Zp^{d'}\cdot\Mcen\Leftrightarrow  \nonumber\\
\forall j:\;\Zp^{d'}\cdot\Mab C_j \alpha & \subseteq & \Zp^{d'}\cdot\mbox{diag}(p^{r_1},p^{r_2},\dots,p^{r_{d'-1}},1)\Leftrightarrow  \nonumber\\
\forall j,\;\forall i \in \{1,\dots,d'-1\}:\;  \Mab C_j\alpha^i &\equiv& 0 \mbox{ mod } p^{r_i} \Leftrightarrow \nonumber \\
\forall i: \; \Mab M(\alpha^i) &\equiv& 0 \mbox{ mod } p^{r_i}. \label{newlatticecondition}
\end{eqnarray}
To be admissible for~$\Lcen$, the lattice~$\Lab$ has to be contained in the solution space of a system of linear equations in~$\mathbb{Z}/(p^{r_1})$. Here we interpret congruences mod $p^{r_i}$ as congruences mod $p^{r_1}$ in the obvious way: $x\equiv0\mod p^{r_i}\Leftrightarrow p^{r_1-r_i}x\equiv0\mod p^{r_1}$. Let~$\{e_j\}_{1\leq j\leq d}$ denote the elementary divisors of this system of linear equations. Then -~with respect to suitable coordinates for the abelianisation~$\LieL_p/{\LieL_p'}$~- condition (\ref{newlatticecondition}) reads as
\begin{equation*}
\Mab\,\mbox{diag}(p^{e_1},\dots,p^{e_d})\equiv 0 \mbox{ mod }p^{r_1}.
\end{equation*}
As $\phi_{\LieL}$ is full, $e_i\leq r_1$ for all $i$ and thus 
\begin{equation*}
|\LieL_p:X(\Lcen)|=p^{\sum_{i=1}^{d}(r_1-e_i)}. \label{formula-for-index-of-X(Lcen)}
\end{equation*}
 
If, in general, $\Lcen=p^r\Lcen^{max}$, where~$\Lcen^{max}$ is the maximal element in the homothety class, the same analysis applies if we replace $r_i$ by $r_i+r$ in~(\ref{newlatticecondition}). It is thus easy to see that 
\begin{equation}
|\LieL_p:X(\Lcen)|=p^{rd}|\LieL_p:X(\Lcen^{max})|.\label{formula:hom. class}
\end{equation}

\begin{definition} \label{definition:new-weight-function}
Let $[\Lcen]$ be a vertex in $\D_{d'}$. Set
\begin{eqnarray*}
 w([\Lcen])&:=&log_p(|\LieL_p':\Lcen^{max}|),\label{weight_function}\\
w'([\Lcen])&:=& w([\Lcen])+\mbox{log}_p(|\LieL_p:X(\Lcen^{max})|).\label{new_weight_function}
\end{eqnarray*}
where $\Lcen^{max}$ denotes the $\leq$-maximal element in $[\Lcen]$.
\end{definition}

We can now express $\z_{\LieL,p}(s)$ via a generating function associated to a weight function on the vertices of the building $\Delta_{d'}$.
\begin{lemma}\label{lemma 1} Let  
\begin{equation*}
A(p,T):=\sum_{[\Lcen]}p^{w([\Lcen])\cdot d}T^{w'([\Lcen])} \label{new_generating_function}
\end{equation*}
where~$\Lcen$ ranges over the full lattices within the derived ring~$\LieL'$. Then

\begin{equation*}
\boxed{\zeta_{\LieL_p}^\triangleleft(s)=\zeta_{\Zp^d}(s)\zeta_p((d+d')s-dd')A(p,p^{-s}).\label{reduction to weight function}}
\end{equation*}
\end{lemma}
\begin{proof} Combine (\ref{GSS-Lemma6.1}), (\ref{formula:hom. class}) and Definition \ref{definition:new-weight-function}.
\end{proof}
We are thus left with the study of the generating function $A(p,T)$. For future reference we note the following
\begin{observation}\label{observation:fu.eq.}
\begin{eqnarray}
\z_{\LieL,p}(s)|_{p\rightarrow p^{-1}}&=&(-1)^{d+d'}p^{\binom{d+d'}{2}-(2d+d')s}\cdot\z_{\LieL,p}(s) \Leftrightarrow \label{fun.eq.}\\
A(p,T)|_{\substack{p \rightarrow p^{-1} \\T \rightarrow T^{-1}}}&=&(-1)^{d'-1}p^{\binom{d'}{2}}\cdot A(p,T). \label{fun.eq.buildings}
\end{eqnarray}
\end{observation}

\section{\T-groups of class $2$ with centres of rank $2$}\label{section3}

In this section we will prove Theorem~\ref{general rank 2 centres theorem}. To that purpose we will recall in~\ref{classification} the classification of class~$2$ \T-groups with centres of Hirsch length $2$ up to commensurability given by Grunewald and Segal. In~\ref{building blocks} we compute the normal zeta functions of the ``building blocks'' in this classification. This will leave us well prepared to prove Theorem~\ref{general rank 2 centres theorem} in~\ref{proof of theorem 2}.

\subsection{Classification up to commensurability}\label{classification}
In \cite{GSegal/84} Grunewald and Segal give a classification of {\sl radicable}~\T-groups of class~$2$ and of finite Hirsch length with centres of Hirsch length $2$, called $\mathcal{D}^*$-groups. As this classification is just up to commensurability we may assume that~$Z(G)/G'$ is {\sl free} abelian of rank~$\leq 1$, where~$Z(G)$ is the centre of the group. As the well known classification of alternating bilinear forms yields the classification of \T-groups with~$G'$ cyclic (cf.~\cite{GScharlau/79}) one may assume that this quotient is indeed trivial. In other words we may assume that the alternating bilinear map

\begin{eqnarray}
\phi_{G}:G/Z(G)\times G/Z(G)&\rightarrow&Z(G) \label{bilinear map - groups} \\
(aZ(G),bZ(G))&\mapsto&[a,b] \nonumber
\end{eqnarray}
is full.

To make this paper self-contained we recall a  Definition and Theorems~$6.2$ and~$6.3$ of Grunewald and Segal~(\cite{GSegal/84}) with slightly adjusted notation. 

\begin{definition}
Let $G$ be a \Dstar-group. A central decomposition of $G$ is a family~$\{H_1,\dots,H_m\}$ of subgroups of $G$ such that
\begin{enumerate} 
\item $Z(H_i)=Z(G)$ for each i;
\item ~$G/Z(G)$ is the direct product of the subgroups $H_i/Z(G)$; and 
\item $[H_i,H_j]=1$ whenever $i\not=j$. The group $G$ is (centrally) indecomposable if the only such decomposition is $\{G\}$.
\end{enumerate}
\end{definition}
\begin{theorem} \label{GS-theorem 1} Every \Dstar-group $G$ has a central decomposition into decomposable constituents, and the decomposition is unique up to an automorphism of~$G$. In particular, the constituents are unique up to isomorphism.
\end{theorem}

\begin{theorem}  \label{GS-theorem 2}
(i) Let $G$ be an indecomposable \Dstar-group of Hirsch length $n+2$. Then, with respect to a suitable basis of $G/Z(G)$ and a suitable basis~$(y_1,y_2)$ of~$Z(G)$, the map~$\phi_G$ is represented by a matrix $M({\bf y})$ as follows:
\begin{itemize}
\item $n=2r+1.$ $$M({\bf y})=M^{r}_0({\bf y})=\left( \begin{array}{cc}
0&{\bf B} \\
-{\bf B}^t&0
\end{array}
\right),$$ where
\begin{eqnarray}
{\bf B}={\bf B}({\bf y})&=&\left( \begin{array}{ccccc}
                                  y_2 & 0 & 0&\ldots & 0 \\
                                  y_1 & y_2 & 0&\ldots & 0 \\
                                  0 & y_1 & y_2&\ldots & 0 \\
                                  0 & 0 & y_1&\ldots & 0 \\
                                  0  &   & & \vdots  &   \\
                                  0 & 0 & 0&\ldots & y_2 \\
                                  0 & 0 & 0&\ldots & y_1 
      \end{array} \right)_{(r+1)\times r}\label{indecomposable-odd-presentation}
\end{eqnarray}

\item $n=2r.$ $$M({\bf y})=M_{(f,e)}({\bf y})=\left( \begin{array}{cc}
0&{\bf B} \\
-{\bf B}^t&0
\end{array}
\right),$$ where
 \begin{eqnarray}
{\bf B}={\bf B}({\bf y})&=&\left( \begin{array}{cccccc}
                         y_1+a_1y_2 & y_2 & 0 & 0&\ldots & 0 \\
                            -a_2y_2 & y_1 & y_2 & 0&\ldots & 0 \\
                             a_3y_2 & 0 & y_1 & y_2&\ldots & 0 \\
                            -a_4y_2 & 0 & 0 & y_1&\ldots & 0 \\
                             \vdots  &   &   & & \vdots  &   \\
                   (-1)^ra_{r-1}y_2 & 0 & 0 & 0&\ldots & y_2 \\
                   (-1)^{r+1}a_{r}y_2 & 0 & 0 & 0&\ldots & y_1 
      \end{array} \right)_{r\times r}\label{indecomposable-even-presentation}
\end{eqnarray}
and $det({\bf B}({\bf y}))=g(y_1,y_2)=y_1^r+a_1y_1^{r-1}y_2+\dots+a_ry_2^r\in\mathbb{Q}[y_1,y_2]$ \label{polynomial_equation} is such that $g(y_1,1)$ is primary, say $g=f^e$ for $f$ irreducible over $\mathbb{Q}$, $e\in \mathbb{N}$.
\end{itemize}
(ii) If $G$ is any \Dstar-group, then with respect to a suitable basis as above, $\phi_G$ is represented by the diagonal sum of matrices like $M({\bf y})$ above.
\end{theorem}

\subsection{Zeta functions of indecomposable \Dstar-groups}\label{building blocks}
We will now compute the local normal zeta functions of indecomposable \Dstar-groups of odd and even Hirsch length in Propositions~\ref{indecomposable_odd} and~\ref{indecomposable_even}, respectively.

\begin{proposition} \label{indecomposable_odd}
Assume $G$ is an indecomposable \Dstar-group of Hirsch length $2r+1+2$ as in the statement of Theorem~\ref{GS-theorem 2} above with $r\geq 1$. Then for all primes $p$
\begin{equation}
A(p,T)=\frac{1+p^{2r+1}T^{2r+1}}{1-p^{2r+2}T^{2r+1}}. \label
{formula indecomposable odd} 
\end{equation}
\end{proposition}

\begin{proof} Recall from Lemma~\ref{lemma 1} that 
$$A(p,T)=\sum_{[\Lcen]}p^{w([\Lcen])\cdot(2r+1)}T^{w'([\Lcen])}$$
for weight functions $w,w'$ defined in Definition~\ref{definition:new-weight-function} on vertices~$[\Lcen]$ of the Bruhat-Tits tree. The integer $w([\Lcen])$ is just the vertices' distance from the root vertex, i.e. the length of the shortest path from~$[\Zp^2]$ to~$[\Lcen]$.
For a given maximal lattice $\Lcen$ of elementary divisor type $(p^{s},1)$, $s \geq 0$, there is an element $\alpha =(\alpha_i^j)\in Sl_2(\Zp)$ such that the admissibility condition (\ref{newlatticecondition}) becomes 
\begin{equation*}
\Lab M_0^r(\alpha^1 )\equiv 0 \mbox{ mod } p^s 
\end{equation*}
Notice that the $(2r+1)\times(2r+1)$-matrix $M_0^r(\alpha^1 )$ has always determinant zero and a $2r$-minor which is a $p$-adic unit which follows from inspection of the matrix in~(\ref{indecomposable-odd-presentation}). Recall that by Definition~\ref{definition:new-weight-function} the integer~$w'([\Lcen])$ measures the index of $\Lcen$ in $\LieL_p'$ and the largest lattice in the abelian quotient admissible for $\Lcen$. In the present case it is easy to verify that
$$\boxed{w'([\Lcen])=s+2rs.}$$ 
The $p+1$ sector-families in the Bruhat-Tits tree are just the sub-trees generated by the root vertex together with one of its neighbours and all of its descendants. They overlap in their boundary, the root vertex. The generating function restricted to the {\sl interior} of any of the~$p+1$ sector-families will therefore equal 
$$\frac{p^{2r+1}T^{2r+1}}{1-p^{2r+2}T^{2r+1}}.$$ Summing over all $p+1$ sector-families gives 
\begin{equation*}
A(p,T)=1+(p+1)\frac{p^{2r+1}T^{2r+1}}{1-p^{2r+2}T^{2r+1}}=\frac{1+p^{2r+1}T^{2r+1}}{1-p^{2r+2}T^{2r+1}}
\end{equation*} as claimed. 
\end{proof}

We apply Observation~\ref{observation:fu.eq.} to equation~(\ref{formula indecomposable odd}) to deduce
\begin{corollary}\label{fun.eq. for indecomposable_odd}
For all primes $p$
$$\z_{G,p}(s)|_{p\rightarrow p^{-1}}=-p^{\binom{2r+3}{2}-(4r+4)s}\cdot\z_{G,p}(s)$$
\hfill$\Box$
\end{corollary}

Indecomposable \Dstar-groups of even Hirsch length are more interesting: 

\begin{proposition} \label{indecomposable_even}
Assume $G$ is an indecomposable \Dstar-group of Hirsch length $2r+2$ as in Theorem~\ref{GS-theorem 2} above with $g=f^e$, $f$ irreducible over~$\Q$, $e\in \mathbb{N}$, and let $p$ be a prime unramified in $\Q[t]/f(t)$. Let $\fFp$ be the number of distinct {\sl linear} factors in $\overline{f(t)}$, the reduction of $f(t)$ mod $p$. Then 
\begin{equation*}
A(p,T)=\frac{P_1(p,T)+\fFp P_2(p,T)}{(1-p^{2r+1}T^{2r-1})(1-p^{2r+1}T^{2r+1})(1-p^{(2r+1)e-1}T^{(2r-1)e})},\label{formula:indecomposable_even}
\end{equation*}

where
\begin{eqnarray}
P_1(p,T)&=&(1-p^{2r+1}T^{2r-1})(1+p^{2r}T^{2r+1})(1-p^{(2r+1)e-1}T^{(2r-1)e}), \nonumber \\
P_2(p,T)&=&p^{2r}T^{2r-1}(1-T^2)(1-p^{(2r+1)e}T^{(2r-1)e})\label{formula: W_i indecomposable even}.
\end{eqnarray}
\end{proposition}

\begin{remark} It is well-known that, for almost all primes $p$, the reduction mod $p$ of the minimal polynomial of a primitive element of a finite extension~$K/\mathbb{Q}$ determines the decomposition behaviour of the prime $p$ in the ring of integers~$\oK$ (cf \cite{Neukirch/92}, Satz (8.3)). Du Sautoy and Grunewald~(\cite{duSG/00}), however, proved the existence of associated varieties whose reduction behaviour mod~$p$ governs the local zeta functions. The two questions ``overlap'' if these varieties are zero-dimensional.
\end{remark}

\begin{proof} We note that $d=2r=2e\deg(f)$, $d'=2$. For a given maximal lattice $\Lcen$ of type $(p^{s},1)$, $s \geq 1$, there exists $\alpha =(\alpha_i^j)\in Sl_2(\Zp)$ such that that the admissibility condition (\ref{newlatticecondition}) reads as 
\begin{equation*}
\Lab M_{(f,e)}(\alpha^1 )\equiv 0 \mbox{ mod } p^s. 
\end{equation*}
The column vector $\alpha^1=(\alpha_1^1,\alpha_2^1)^t$ is unique only up to multiplication by a $p$-adic unit and up to addition of multiples of $p^s\cdot\alpha^2$ and therefore determines a unique element in $\mathbb{P}^1(\mathbb{Z}/(p^s)).$ In particular, $v_p(f(\alpha^1)) \leq s,$ where $v_p$ is the $p$-adic valuation. 
We say that $$[\Lcen]\equiv x\;:\Leftrightarrow \alpha^1 \equiv x \mbox{ mod }p, \;x\in \mathbb{P}^1(\mathbb{F}_p).$$ 
We have
\begin{eqnarray}
 A(p,T)&=&1+\sum_{x\in \mathbb{P}^1(\mathbb{F}_p)} \sum_{[\Lcen]\equiv x}p^{w([\Lcen])\cdot2r}T^{w'([\Lcen])}\nonumber \\
       &=&\sum_{x\in \mathbb{P}^1(\mathbb{F}_p)}\underbrace{\left(\frac{1}{p+1}y+\sum_{[\Lcen]\equiv x}p^{w([\Lcen])\cdot2r}T^{w'([\Lcen])}\right)}_{=:A(p,T,x)}.\nonumber
\end{eqnarray}
Here we have distributed the constant term of $A(p,T)$ over the~$p+1$ sector-families to get a formal sum over the points in $\mathbb{P}^1(\mathbb{F}_p)$.

For the computation of the $A(p,T,x)$ we only have to distinguish two cases:
\begin{description}
\item[Case 1:]$ f(x)\not\equiv 0 \mbox{ mod } p.$
\item[Case 2:]$ f(x)\equiv 0 \mbox{ mod } p,$ i.e. $x$ is a {\sl simple\footnote{The zero has to be {\sl simple} since we have excluded ramified primes.} zero} mod $p$.

\end{description}

As to case 1: $M_{(f,e)}(\alpha^1 )$ has unit determinant in this case, so by Definition~\ref{definition:new-weight-function}
$$w'([\Lcen])=s(2r+1).$$ We readily compute 
\begin{equation*}
A(p,T,x)=\frac{1}{p+1}+\frac{p^{2r}T^{2r+1}}{1-p^{2r+1}T^{2r+1}}=:A_{\emptyset}(p,T) \label{A_emptyset}.
\end{equation*}

As to case 2: By Hensel's Lemma we know that a simple zero of $f(t)$ mod~$p$ lifts to a zero mod $p^n$ for all $n$ and ``eventually'' to a zero in $\Zp$. Vertices of the Bruhat-Tits tree with distance $n$ from the root vertex are identified with points on the projective line $\mathbb{P}^1(\mathbb{Z}/p^n)$. A zero in~$\mathbb{P}^1(\Qp)$ may hence be viewed as an end in the Bruhat-Tits tree.

 We know that $f$ has a simple zero mod $p$ at $x\equiv\alpha^1$. However, $w'([\Lcen])$ depends on the exact vanishing order of $f(\alpha^1)$ mod $p^s$. 
By Definition~\ref{definition:new-weight-function} we have
\begin{equation}
\boxed{w'([\Lcen])=s+2(rs-\mbox{min}\{s,e\cdot v_p(f(\alpha^1))\}).} \label{weight irreducible, non-unit case}
\end{equation}
This follows from the fact that the matrix $M(\alpha^1)$ always has a $2(r-1)\mbox{-minor}$ which is a $p$-adic unit, which again follows from inspection of the presentation~(\ref{indecomposable-even-presentation}).

Thus the map
$$[\Lcen] \mapsto p^{w([\Lcen])\cdot2r}T^{w'([\Lcen])} $$
factorizes over the set $\coneN:=\{(a,b)\in \mathbb{N}_{>0}^2|\; a\geq b \geq 1\}$ as $\psi \phi$ where
\begin{eqnarray}
\phi:[\Lcen] &\mapsto & (s, v_p(f(\alpha^1))) \nonumber \\
\psi:(a,b) &\mapsto &p^{2ra}T^{(2r+1)a-2min\{a,eb\}}. \label{formula-factorization-over-cone}
\end{eqnarray}
We view $N$ as intersection of~$\N_{>0}^2$ with a closed rational polyhedral cone~$C$ in~$\mathbb{R}_{>0}^2$.
Notice that the cardinality of the preimage of $\phi$ is given by
\begin{equation} \label{preimage phi} 
|\phi^{-1}(a,b)|=\left\{\begin{array}{cl}
                  1 & \mbox{if }a=b, \\
                  (1-p^{-1})p^{a-b} & \mbox{if }a>b.  
                 \end{array} 
                \right.
\end{equation}
What we want to compute is 
\begin{equation*}
A(p,T,x)=\frac{1}{p+1}+\sum_{(a,b) \in N} |\phi^{-1}(a,b)| \psi(a,b). \label{gen.fun.over S}
\end{equation*}
In order to eliminate the ``$\min$'' in the expression~(\ref{formula-factorization-over-cone}) we decompose the sector~$N$ into sub-cones $N_j$ on which it is easier to sum over~$|\phi^{-1}(a,b)| \psi(a,b)$. The technical idea is to substitute the variables $X,Y$ in the zeta functions for the respective sub-cones $N_j$ by certain Laurent monomials\footnote{By a {\sl Laurent monomial} in $p,T$ we mean of course a term $p^aT^b$ where $a,b\in\Z$.} $m_{jX}(p,T)$, $m_{jY}(p,T)$ in $p$ and $T$.
We choose the decomposition
\begin{eqnarray*}
N&=&N_0+N_1 + N_2, \\
N_0&:=& \{(a,b)\in N|\; a=b\geq 1\},\nonumber \\
N_1&:=& \{(a,b)\in N|\; eb>a>b\geq1\}, \nonumber \\
N_2&:=& \{(a,b)\in N|\; a \geq eb\}, \nonumber
\end{eqnarray*}
with respective zeta functions
\begin{eqnarray}
F_0(X,Y)&:=&\sum_{\substack{(a,b)\in N_0}}X^aY^b=\frac{YX}{1-YX},\nonumber\\
F_1(X,Y)&:=&\sum_{\substack{(a,b)\in N_1}}X^aY^b=\frac{YX^2}{(1-X)(1-YX)}-\frac{YX^e}{(1-X)(1-YX^e)}, \nonumber \\
F_2(X,Y)&:=&\sum_{\substack{(a,b)\in N_2}}X^aY^b=\frac{YX^e}{(1-X)(1-YX^e)}. \nonumber
\end{eqnarray}
We set $n_j:=\dim(N_j)$, $j=0,1,2$, where, of course, by the {\sl dimension} of~$N_j$ we mean the dimension of the corresponding polyhedral cone in~$\mathbb{R}_{>0}^2$. The Laurent monomials in $p,T$ which we have to substitute for $X,Y$ are easily read off from formulae~(\ref{formula-factorization-over-cone}) and~(\ref{preimage phi}). We have
\begin{eqnarray}
\lefteqn{A(p,T,x)=\frac{1}{p+1}+\sum_{j=0}^2(1-p^{-1})^{n_j-1}F_j(X,Y)|_{\substack{X=m_{jX}(p,T)\\Y=m_{jY}(p,T)}}}\nonumber\\
  &=&\frac{1}{p+1}+F_0(p^{2r}T^{2r+1},T^{-2}) \nonumber \\
         &&\quad +(1-p^{-1})\left(F_1(p^{2r+1}T^{2r-1},p^{-1})+F_2(p^{2r+1}T^{2r+1},p^{-1}T^{-2e})\right) \nonumber\\
         &=:&A_{\{1\}}(p,T).  \label{A_1}
\end{eqnarray}

The formula for $A(p,T)$ is now obtained by counting occurrences of the two cases and routine computations with rational functions. Indeed, we have 
\begin{eqnarray}
\lefteqn{A(p,T)=\left(p+1-\fFp \right)A_{\emptyset}(p,T)+\fFp A_{\{1\}}(p,T)=}\nonumber \\
&&\frac{P_1(p,T)+\fFp P_2(p,T)}{(1-p^{2r+1}T^{2r-1})(1-p^{2r+1}T^{2r+1})(1-p^{(2r+1)e-1}T^{(2r-1)e})},\nonumber 
\end{eqnarray}
with $P_i(p,T)$, $i=1,2$ defined as in the statement of Proposition~\ref{indecomposable_even}.
\end{proof} 

We apply Observation~\ref{observation:fu.eq.} to equation~(\ref{formula: W_i indecomposable even}) to deduce
\begin{corollary}\label{fun.eq. for indecomposable_even}
For almost all primes $p$
$$\z_{G,p}(s)|_{p\rightarrow p^{-1}}=-p^{\binom{2r+2}{2}-(4r+2)s}\cdot\z_{G,p}(s)$$
\end{corollary}
\hfill$\Box$

\subsection{Proof of Theorem~\ref{general rank 2 centres theorem}}\label{proof of theorem 2}

As commensurable \T-groups have identical local zeta functions for almost all primes we may assume $G$ is a \Dstar-group and that $\LieL(G)$ is presented as in (\ref{Lie ring presentation}) where 

\begin{eqnarray}
M({\bf y})&:=&\bigoplus_{i=1}^m M_{(f_i,{\bf e}_i)}({\bf y})\oplus \bigoplus_{k=1}^n \underbrace{M_0^{l_k}({\bf y})}_{\substack{\mbox{``indecomp. of odd}\\ \mbox{Hirsch length''} }}, \label{general Dstar Lie algebra presentation formula} \\
M_{(f_i,{\bf e}_i)}({\bf y})&:=&\oplus_{j=1}^{r_i} \underbrace{M_{(f_i,e_{ij})}({\bf y})}_{\substack{\mbox{``indecomp. of even}\\ \mbox{Hirsch length''} }},\nonumber
\end{eqnarray}
for a set of irreducible polynomials $\mathcal{F}=\{f_1,\dots,f_m\}$ of degree~$deg(f_i)=d_i$ and suitable $m,n,l_k,e_{ij},\,r_i\in\N$, using the notation of Theorem~\ref{GS-theorem 2}. The vector~${\bf e}_i=(e_{ij})$ will be called the {\rm vector of multiplicities of $f_i$}.

Note that~(\ref{general Dstar Lie algebra presentation formula}) just spells out $(ii)$ of Theorem~\ref{GS-theorem 2} and that 
$$d=\sum_{i=1}^{m}2d_i\left(\sum_{j=1}^{r_i}e_{ij}\right)+\sum_{k=1}^n(2l_k+1),\: d'=2.$$ 

Excluding a finite set of primes if necessary we may assume that $p$
is unramified in each of the $\Q[t]/f_i(t)$, $f_i\in\mathcal{F},$ and
that if any two polynomials $f_i,f_j\in\mathcal{F}$ have simple zeros
$\alpha_i$ and $\alpha_j$ in $\Zp$ and $\alpha_i\equiv\alpha_j\mod p$,
then $\alpha_i=\alpha_j$. Viewing roots as ends in the Bruhat-Tits
tree emanating from the root gives an easy geometric interpretation  this condition: distinct ends shall only have the root in common. 

Again we have to analyse
\begin{equation*}
A(p,T)=\sum_{x\in \mathbb{P}^1(\mathbb{F}_p)}\underbrace{\left(\frac{1}{p+1}+\sum_{[\Lcen]\equiv x}p^{w([\Lcen])\cdot2r}T^{w'([\Lcen])}\right)}_{=A(p,T,x)}.
\end{equation*}

In fact we will see that, for~$x\in \mathbb{P}^1(\mathbb{F}_p)$, $A(p,T,x)$ is a rational function which only depends on the set of~$f_i$'s, together with their multiplicities, which happen to have a simple zero at~$x$, and on $\sum_{k=1}^nl_k$ and~$n$, the total size and number of ``indecomposable blocks of odd Hirsch length''. In particular the dependence of $A(p,T,x)$ {\sl on the prime $p$} will be encoded in the set $$I(x):=\{i\in \{1,\dots,m\}|\,f_i(x)\equiv 0 \mbox{ mod }p \}$$ so that for~$I\in\{1,\dots,m\}$ we may set
\begin{equation}
A_{I}(p,T):=\left\{ \begin{array}{ll}
                     A(p,T,x) & \mbox{ if } I=I(x)\mbox{ for some }x\in\mathbb{P}^1(\Fp),\\
                     0 &\mbox{ otherwise.}
                   \end{array} \right.\label{A(p,T,x) depends only on I(x)}
\end{equation} 
Assuming for a moment that~(\ref{A(p,T,x) depends only on I(x)}) is well-defined we have 
\begin{eqnarray}
A(p,T)&=&\sum_{x\in \mathbb{P}^1(\mathbb{F}_p)} A(p,T,x)=\nonumber\\
       &&\sum_{I\subseteq\{1,\dots,m\}}c_{p,I}A_{I}(p,T)=\label{c_p,I-formula}\\
       &&(p+1)A_{\emptyset}(p,T)+\sum_{I\subseteq\{1,\dots,m\}}c_{p,I}\left(A_{I}(p,T)-A_{\emptyset}(p,T)\right), \label{functional-equation-formula}
\end{eqnarray}
where $c_{p,I}$ are defined as in Theorem~\ref{general rank 2 centres theorem}. Note that $\sum_{I\subseteq\{1,\dots,m\}}c_{p,I}=p+1$ trivially. Together with Lemma~\ref{lemma 1}, formula~(\ref{c_p,I-formula}) clearly implies the first part of Theorem~\ref{general rank 2 centres theorem}, whereas the functional equation will follow from~(\ref{functional-equation-formula}) and Observation~\ref{observation:fu.eq.} once we have established that
\begin{equation*}
A'_{I}(p,T):=\left\{\begin{array}{rl}
                     (p+1)A_{\emptyset}(p,T) & \mbox{ if }\emptyset=I ,y\\
        A_{I}(p,T)-A_{\emptyset}(p,T)& \mbox{ if }\emptyset\not=I\subseteq \{1,\dots,m\} 
                     \end{array}\right.
\end{equation*}
satisfies the functional equation
\begin{equation}
A'_I(p,T)|_{\substack{p \rightarrow p^{-1} \\T \rightarrow T^{-1} }}=-p \cdot A'_{I}(p,T),\; I\subseteq\{1,\dots,m\}. \label{functional-equation_A'_I(p,T)}
\end{equation}

We now proceed to compute $A_{\emptyset}(p,T)$ and show that our definition~(\ref{A(p,T,x) depends only on I(x)}) is well-defined for a general $\emptyset\not=I\subseteq\{1,\dots,m\}.$ We will give an algorithm for the computation of the~$A_{I}(p,T)$, $I\subseteq\{1,\dots,m\}$ and prove inductively that the rational functions~$A'_{I}(p,T)$ satisfy the functional equation~(\ref{functional-equation_A'_I(p,T)}).\bigskip

To compute $A(p,T,x)$ for $x$ such that~$f_i(x)\not\equiv 0$ mod~$p$ for all~$1\leq i \leq m$ fix a class~$[\Lcen]\equiv x$ of distance~$s$ from the root. Consider the admissibility condition
\begin{equation*}
\Lab M(\alpha^1 )\equiv 0 \mbox{ mod } p^s. 
\end{equation*}
By our choice of $x$ all the $M_{(f_i,{\bf e}_i)}(\alpha^1)$, $1\leq k \leq n$, have unit determinant. All the $M_0^{l_k}(\alpha^1)$ have determinant zero but have a unit minor of maximal size. It follows that 
$$w'([\Lcen])=s(d+1-n)$$
and therefore
\begin{eqnarray*}
A'_{\emptyset}(p,T)&=&(p+1)A_{\emptyset}(p,T)=\frac{1+p^dT^{d+1-n}}{1-p^{d+1}T^{d+1-n}}.
\end{eqnarray*}
Thus, for $I=\emptyset$, one easily reads off the functional equation~(\ref{functional-equation_A'_I(p,T)}). \bigskip

Now we consider $A(p,T,x)$ for an $x$ such that 
$$\exists\;1\leq i \leq m:f_i(x)\equiv 0 \mbox{ mod }p$$
 and we set $I(x):=\{i\in \{1,\dots,m\}|\,f_i(x)\equiv 0 \mbox{ mod }p\}\not=\emptyset$. Fix a class $[\Lcen]\equiv x$ of distance $s$ from the root. It is determined by an element~$\alpha^1\in\mathbb{P}^1(\mathbb{Z}/(p^s))$ congruent to $x$ mod $p$. We may assume that locally around~$x$ all the $f_i(t)$, $i\in I(x)$, are given by $u_it=0$ for some~$p$-adic units~$u_i\in\Zp^*$. We put~$v_p(\alpha^1)=c$, $1\leq c \leq s$. Then 
\begin{equation*}
w'([\Lcen])=s(d+1-n)-2(\sum_{i\in I(x)}\sum_{j=1}^{r_i}\min(s,e_{ij}\cdot c)).
\end{equation*}   
(We recover formula (\ref{weight irreducible, non-unit case}) as a special case.)
Exactly as in the proof of Theorem~\ref{indecomposable_even}, the map
$$[\Lcen] \mapsto p^{w([\Lcen])\cdot d}T^{w'([\Lcen])} $$
factorizes over the set $N:=\{(a,b)\in \mathbb{N}_{>0}^2|\; a\geq b \geq 1\}$ as $\psi \phi$ where
\begin{eqnarray}
\phi:[\Lcen] &\mapsto & (s,c) \nonumber \\
\psi:(a,b) &\mapsto &p^{ad}T^{a(d+1-n)-2(\sum_{i\in I}\sum_{j=1}^{r_i}\min(a,e_{ij}b))}. \nonumber
\end{eqnarray}
We have 
\begin{equation*}
A(p,T,x)=\frac{1}{p+1}+\sum_{(a,b)\in N} |\phi^{-1}(a,b)| \psi(a,b), \label{gen.fun.over N}
\end{equation*}
a rational function in $p$ and $T$ which depends only on $({\bf e}_i)_{i\in I(x)}$,~$d$ and~$n$. This shows that~$A_I(p,T)$ is well-defined and proves the first part of Theorem~\ref{general rank 2 centres theorem}. 

We will now give an algorithm to compute the~$A_I(p,T)$ and prove the functional equation. To ease notation we will rewrite $\psi$ such that
\begin{equation}
\psi(a,b)=p^{ad}T^{c_0b+\sum_{r=1}^sc_rmin(a,e_rb)+c_{s+1}a} \label{rewrite psi}
\end{equation}
with constants $c_i,e_r\in \mathbb{Z},\; i\in\{0,\dots, s+1\},\; r\in\{1,\dots,s\},\;e_r> 1$. We will apply the same trick as in the proof of Proposition~\ref{indecomposable_even}: To eliminate the ``min'' in~(\ref{rewrite psi}) we again define sub-sectors $N_j$, $0 \leq j \leq s+1$, of $N$ on which $|\phi^{-1}(a,b)| \psi(a,b)$ are easier to sum over. We choose the subdivision

\begin{eqnarray}
N&=&N_0+N_1+\dots+N_{s+1}, \nonumber \\
N_0&:=&\{ (a,b)\in N|\;b=a \geq 1\}, \nonumber \\
N_1&:=&\{ (a,b)\in N|\; e_1b >a\geq b\geq 1\}, \nonumber \\
N_2&:=&\{ (a,b)\in N|\; e_2b > a \geq e_1b\}, \nonumber \\
   &\vdots& \nonumber \\
N_{s+1}&:=&\{ (a,b)\in N|\;a \geq e_sb\}. \nonumber
\end{eqnarray}
We put $n_j:=\dim(N_j)$ and observe that~$n_0=1$ and~$n_j=2$ for~$1\leq j \leq s+1$. Let $$F_j(X,Y):=\sum_{(a,b)\in N_j}X^aY^b,\;0\leq j \leq s+1$$ denote the generating function for the sub-sector $N_j$. We can now break up the generating function $A_I(p,T)$ as follows:
\begin{eqnarray*}
A_I(p,T)=A(p,T,x)&=&\frac{1}{p+1}+\sum_{j=0}^{s+1} \sum_{(a,b)\in N_j} |\phi^{-1}(a,b)| \psi(a,b)=\\
        &&\frac{1}{p+1}+\sum_{j=0}^{s+1}(1-p^{-1})^{n_j-1}F_j(X,Y)|_{\substack{X=m_{jX}(p,T)\\Y=m_{jY}(p,T)}},
\end{eqnarray*}
where $m_{jX}(p,T),m_{jY}(p,T)$, $0\leq j \leq s+1$, are suitable Laurent monomials. They are recorded in Table~\ref{table-general-dstar-group}, together with the $n_j$ and the generating functions~$F_j(X,Y)$. 

\begin{table}[h]
\centering
\setlength{\extrarowheight}{9pt}
\begin{tabular}{c|c|c|c|c}
j&$n_j$&$F_j(X,Y)$&$m_{jX}(p,T)$&$m_{jY}(p,T)$ \\ \hline 
0&1&$\frac{YX}{1-YX}$&$p^{d+1}T^{\sum_{r=1}^{s+1}c_r}$&$p^{-1}T^{c_0}$ \\
1&2&$\frac{YX^2}{(1-X)(1-YX)}-\frac{YX^{e_1}}{(1-YX^{e_1})(1-X)}$&$p^{d+1}T^{\sum_{r=1}^{s+1}c_r}$&$p^{-1}T^{c_0}$ \\
2&2&$\frac{Y(X^{e_1}-X^{e_2})}{(1-YX^{e_1})(1-YX^{e_2})(1-X)}$&$p^{d+1}T^{\sum_{r=2}^{s+1}c_r}$&$p^{-1}T^{c_0+c_1e_1}$ \\
\vdots&\vdots&\vdots& \vdots\\
s&2&$\frac{Y(X^{e_{s-1}}-X^{e_s})}{(1-YX^{e_{s-1}})(1-YX^{e_s})(1-X)}$&$p^{d+1}T^{c_s+c_{s+1}}$&$p^{-1}T^{c_0+\sum_{r=1}^{s-1}c_re_r}$ \\
s+1&2&$\frac{YX^{e_s}}{(1-YX^{e_s})(1-X)}$&$p^{d+1}T^{c_{s+1}}$&$p^{-1}T^{c_0+\sum_{r=1}^{s}c_re_r}$ \\
\end{tabular}\label{table-general-dstar-group}
\caption{In order to compute $A_I(p,T)$ one subdivides $N$ into sub-sectors $N_j$. 
}
\end{table}

Note that the table contains all necessary data to compute $A_I(p,T)$ explicitly. Presently we are only interested in proving one feature of the resulting rational function - the functional equation of $A_I(p,T)-A_{\emptyset}(p,T)$. We will prove it by an induction on the number of $2$-dimensional sub-sectors of~$N$. Proposition~\ref{indecomposable_even} serves as induction base. Assume the functional equation holds for $s$ such sub-sectors. We want to show that it holds for
\begin{eqnarray}
\lefteqn{A'_I(p,T)=A_I(p,T)-A_{\emptyset}(p,T)=} \nonumber \\
         &=&\sum_{(a,b)\in N}|\phi^{-1}(a,b)| T^{c_0b+\sum_{r=1}^sc_rmin(a,e_rb)+c_{s+1}a}-\frac{p^dT^{c_{s+1}}}{1-p^{d+1}T^{c_{s+1}}}. \nonumber 
\end{eqnarray}
By the induction hypothesis, we know that the rational function
\begin{eqnarray}
\lefteqn{B'_I(p,T):=B_I(p,T)-B_{\emptyset}(p,T)=} \nonumber \\
         &=&\sum_{(a,b)\in N}|\phi^{-1}(a,b)| T^{c_0b+\sum_{r=1}^{s-1}c_rmin(a,e_rb)+(c_s+c_{s+1})a}-\frac{p^dT^{c_s+c_{s+1}}}{1-p^{d+1}T^{c_s+c_{s+1}}} \nonumber
\end{eqnarray}
satisfies the functional equation. Note that $B'_I(p,T)$ does not need to come from a~\Dstar-group. It is clearly enough to show the functional equation for the difference 
\begin{eqnarray}
C'_I(p,T)&:=&A'_I(p,T)-B'_I(p,T)=\nonumber\\
         &&A_I(p,T)-B_I(p,T)-(A_{\emptyset}(p,T)-B_{\emptyset}(p,T)) \label{definition C'}  
\end{eqnarray}
We have 
\begin{eqnarray}
A_{\emptyset}(p,T)-B_{\emptyset}(p,T)&=&\frac{p^dT^{c_{s+1}}(1-T^{c_s})}{(1-p^{d+1}T^{c_{s+1}})(1-p^{d+1}T^{c_s+c_{s+1}})}
\label{A_emptyset-B_emptyset}
\end{eqnarray}
As to $A_I(p,T)-B_I(p,T)$, note that the two functions differ only on the sector~$N_{s+1}$, which makes the computation of their difference a lot easier. Using the last row of the above table and setting $$Z:=p^{(d+1)e_s-1}T^{\sum_{r=0}^sc_re_r+c_{s+1}e_s}$$ we find
\begin{eqnarray}
\lefteqn{A_I(p,T)-B_I(p,T)=}\nonumber\\
&&(1-p^{-1})F_{s+1}(X,Y)|_{\substack{X=p^{d+1}T^{c_{s+1}} \\Y=p^{-1}T^{c_0+\sum_{r=1}^sc_re_r}}}\nonumber\\
&&-(1-p^{-1})F_{s+1}(X,Y)|_{\substack{X=p^{d+1}T^{c_s+c_{s+1}} \\Y=p^{-1}T^{c_0+\sum_{r=1}^{s-1}c_re_r}}}=\nonumber\\
&&\frac{(p-1)Z\cdot p^{d}T^{c_{s+1}}(1-T^{c_s})}{(1-Z)(1-p^{d+1}T^{c_{s+1}})(1-p^{d+1}T^{c_s+c_{s+1}})} \label{A_I-B_I}
\end{eqnarray}
Substituting~(\ref{A_emptyset-B_emptyset}) and~(\ref{A_I-B_I}) in equation~(\ref{definition C'}) we get
\begin{equation*}
C'_I(p,T)=\left(\frac{(p-1)Z}{1-Z}-1\right)\frac{p^dT^{c_{s+1}}(1-T^{c_s})}{(1-p^{d+1}T^{c_{s+1}})(1-p^{d+1}T^{c_s+c_{s+1}})}
\end{equation*}
which clearly satisfies $$C'_I(p,T)|_{\substack{p\rightarrow p^{-1}\\T\rightarrow T^{-1}}}=-p^{-1}\cdot C'_I(p,T).$$
This prove Theorem~\ref{general rank 2 centres theorem}.


\section{Proof of Theorem~\ref{main theorem}}\label{section4}

We will compute the generating function 
$$A(p,T)=\sum_{[\Lcen]}p^{w([\Lcen])\cdot 2r}T^{w'([\Lcen])}$$
associated to the vertices of the building $\Delta_3$ as explained in Section~\ref{section2}. We will put to use the decomposition of $\Delta_3$ into {\sl sector-families}~$\SF$, $F\in\Fpthree$, the simplicial complex of proper subspaces in the finite vector space $\Fp^3$. This decomposition is not disjoint. Sector-families overlap at their {\sl boundaries}~$\partial\SF$, and care has to be taken not to over-count. We write 
\begin{equation*}
\pthreetwo=p^2+p+1,\;\ptwoone=p+1 
\end{equation*}
for the~$p$-binomial coefficients. We thus have
\begin{equation*}
A(p,T)=\sum_{F \in \mathcal{F}(p,3)}\ApTSF
\end{equation*}
where 
\begin{eqnarray}
\ApTSF&:=&\frac{1}{\pthreetwo \ptwoone}+\frac{1}{\ptwoone}\sum_{[\Zp^3]\not=[\Lcen]\in\partial\mathcal{S}_F}p^{w([\Lcen])\cdot2r}T^{w'([\Lcen])}\nonumber\\
&&+\sum_{[\Lcen]\in\mathcal{S}_F^{\circ}}p^{w([\Lcen])\cdot2r}T^{w'([\Lcen])}.\label{flag-formula}
\end{eqnarray}

Note that the terms coming from vertices in the sector-family's boundary are {\sl weighted} by factors, reflecting the fact that they are shared by more than one sector-family.

We will now explain why we decompose $\Delta_3$ into sector-families $\mathcal{S}_F$. 
For a given lattice $\Lcen$ of type $(p^{s+t},p^t,1)$, $s,t \geq 0$ there is a unique coset $$\alpha G_{\nu} \in Sl_3(\Zp)/G_{\nu},$$ 
where $G_{\nu}:=\mbox{Stab}_{Sl_3(\Zp)}(\Zp^3\cdot\mbox{diag}(p^{s+t},p^t,1))$,
such that admissibility condition (\ref{newlatticecondition}) becomes 

\begin{eqnarray}
\Lambda_{ab}M(\alpha^1) &\equiv&  0\mbox{ mod } p^{s+t}  \label{eqn_alpha1} \\
\Lambda_{ab}M(\alpha^2) &\equiv&  0\mbox{ mod } p^{t},  \label{eqn_alpha2} 
\end{eqnarray}
where we denote by $\a^j$ the~$j$-th column of the matrix~$\a$. We are free to multiply both $\a^1,\a^2$ by $p$-adic units and to add multiples of $p^s\a^2$ to $\a^1$ and multiples of $\a^1$ to $\a^2$. In particular, if $s\geq 1$ and $t \geq1$, the reduction mod~$p$ $\overline{a^1}$ defines a {\sl point} and $\langle \overline{a^1}, \overline{a^2} \rangle$ defines a {\sl line} in $\PtwooverFp$, and $[\Lcen]$ defines a vertex in the {\sl interior} of $\mathcal{S}_F$, where $F=( \overline{a^1},\langle \overline{a^1}, \overline{a^2} \rangle)$ in $\Fpthree$.

As explained in Section~\ref{section2}, we now have to analyse the elementary divisors of the systems of linear equations~(\ref{eqn_alpha1}) and~(\ref{eqn_alpha2}). This only depends on how the flag $F$ meets the {\sl degeneracy locus} of the matrix $M({\bf y})$. As we are looking to prove a result for almost all primes $p$, we may assume $p+1>r$ to make sure that no line in $\mathbb{P}^2(\Fp)$ is contained in $\overline{V}$, where $V=(\det(R({\bf y}))=0)$. It thus suffices to compute $A(p,T,F)$ for the following two cases:
\begin{description}
\item[Case 1:] \label{offpoint-offline} $det(R(\overline{a^1}))\not\equiv 0 \mbox{ mod }p$. $\ApTSF=:\Woffoff$.
\item[Case 2:] \label{smoothpoint-offline} $det(R(\overline{a^1}))\equiv 0 \mbox{ mod }p$. $\ApTSF=:\Wsmptoff$.
\end{description}

Here the subscripts denote the flag's relative position to the curve $V$: its point either defines a point {\sl off} the curve or a {\sl smooth point}, whereas the line is never a component of $V$. 

As to case 1: Both matrices in~(\ref{eqn_alpha1}) and (\ref{eqn_alpha2}) are non-singular, so that admissibility conditions~(\ref{eqn_alpha1}) and~(\ref{eqn_alpha2}) become
$$ \Lab \equiv 0 \mbox{ mod } p^{s+t}.$$
We compute the weight function on vertices over this sector-family as
\begin{equation*}
\boxed{w'([\Lcen])={s+2t+2r(s+t)}={(2r+1)s+(2r+2)t}.}\label{weight function:A_off/off}
\end{equation*}
The boundary of a sector-family in $\Delta_3$ falls into three components: the root vertex, (maximal) lattices of type $(p^s,1,1)$, and lattices of type $(p^t,p^t,1)$. For $s,t\geq 1$ there are $p^{2(s-1)}$, $p^{2(t-1)}$ and $p^{s-1}p^{t-1}p^{s+t-1}=p^{2s+2t-3}$ lattices in~$\SF$ of type $(p^s,1,1)$, $(p^t,p^t,1)$ and $(p^{s+t},p^t,1)$, respectively. If we write $\ApTSF$ as in~(\ref{flag-formula}), we therefore get
\begin{eqnarray*}
\lefteqn{\Woffoff=}\nonumber \\ 
&&\frac{1}{\pthreetwo \ptwoone}+\frac{1}{\ptwoone}\left(\sum_{s\geq 1}p^{(2r+2)s-2}T^{(2r+1)s}+\sum_{t\geq 1}p^{(4r+2)t-2}T^{(2r+2)t}\right)+\nonumber\\
&&\quad\sum_{s,t\geq 1}p^{(2r+2)s+(4r+2)t-3}T^{(2r+1)s+(2r+2)t}=\nonumber\\
&&\frac{1+p^{2r}T^{2r+1}+p^{2r+1}T^{2r+1}+p^{4r}T^{2r+2}+p^{4r+1}T^{2r+2}+p^{6r+1}T^{4r+3}}{\pthreetwo\ptwoone(1-p^{4r+2}T^{2r+2})(1-p^{2r+2}T^{2r+1})}.\label{A_off/off}
\end{eqnarray*}

As to case 2: In this case, the matrix $M(\alpha^2)$ in~(\ref{eqn_alpha2}) may be assumed to be non-singular, as the column vector~$\alpha^2$ may be moved along the line $\langle \alpha^1,\alpha^2 \rangle$ and we assume that this line contains a point for which the matrix is not degenerate. We may choose affine local coordinates $(x,y,1)$ around $\alpha^1\in \mathbb{P}^2(\Zp)$ such that~(\ref{eqn_alpha1}), (\ref{eqn_alpha2}) may be written as 

\begin{eqnarray*}
\Lambda_{ab}\left(\begin{array}{cc} 
                    0 &\mbox{diag}(x,1,\dots,1) \\
           -\mbox{diag}(x,1,\dots,1) &0 
                     \end{array} \right)
 &\equiv& 0 \mbox{ mod } p^{s+t} \\
\Lambda_{ab}&\equiv&  0\mbox{ mod } p^{t}
\end{eqnarray*}
($x\in p\Z/(p^s)$). Therefore the weight function~$w'$ is given by 
\begin{equation*}
\boxed{w'([\Lcen])={(2r+1)s+(2r+2)t-2\min(s,v_p(x))}.} \label{weight function:A_smpt/off}
\end{equation*}

We begin by summing over the component of the boundary consisting of lattices of type $(p^s,1,1)$, $s\geq1$. We will do this generalizing the method which we developed to prove Theorem~\ref{general rank 2 centres theorem}. We observe that the map 
$$[\Lcen] \mapsto p^{w([\Lcen])\cdot2r}T^{w'([\Lcen])}$$ 
factorizes over the set $N:=\{(a,b,c) \in \mathbb{N}_{>0}^3|\;a\geq b, a \geq c\}$ -~which we view again as the intersection of~$\mathbb{N}_{>0}^3$ with some closed polyhedral cone~$C$ in~$\mathbb{R}_{>0}^3$~- as $\psi \phi$ where 
\begin{eqnarray}
\phi:[\Lcen] &\mapsto& (s,v_p(x),v_p(y)) \nonumber \\
\psi:(a,b,c) &\mapsto& p^{2ra}T^{(2r+1)a-2\min(a,b)} \label{psi - smpt}
\end{eqnarray}
Note that 
\begin{equation} \label{preimage phi - smpt} 
|\phi^{-1}(a,b,c)|=\left\{\begin{array}{ll}
                  1 & \mbox{if }a=b=c, \\
                  (1-p^{-1})\;p^{a-b} & \mbox{if }a>b, a=c, \\
                  (1-p^{-1})\;p^{a-c} & \mbox{if }a>c, a=b, \\
                  (1-p^{-1})^2p^{2a-b-c} & \mbox{if }a>b, a>c. 
                \end{array} 
                \right.
\end{equation}

\noindent Again we will decompose the cone $N$ into sub-cones~$N_j$ on which the values $|\phi^{-1}(a,b,c)|\psi(a,b,c)$ are easier to sum over.
We choose the decomposition
\begin{eqnarray}
N&=&N_0+N_1+N_2+N_3, \label{subdivision N dim 3 into 4 parts}\\
N_0&:=&\{(a,b,c)\in N|\;a=b=c\geq 1\},\nonumber\\
N_1&:=&\{(a,b,c)\in N|\;a=c>b\geq 1\},\nonumber \\
N_2&:=&\{(a,b,c)\in N|\;a=b>c\geq 1\},\nonumber \\
N_3&:=&\{(a,b,c)\in N|\;a>b\geq 1,\,a>c\geq 1\}. \nonumber 
\end{eqnarray}
Again we set $n_j:=\dim(N_j)$. Table~\ref{smpt-table} records the generating functions $F_j(X,Y,Z)$ together with the integers~$n_j$ and Laurent monomials~$m_{jX}(p,T)$, $m_{jY}(p,T)$, $m_{jZ}(p,T)$. The latter are easily read off from~(\ref{psi - smpt}) and~(\ref{preimage phi - smpt}).



\begin{table}[h]
\centering
\setlength{\extrarowheight}{9pt}
\begin{tabular}{c|c|c|c|c|c}
j&$n_j$&$F_j(X,Y,Z)$&$m_{jX}(p,T)$&$m_{jY}(p,T)$&$m_{jZ}(p,T)$ \\ \hline 
0&1&$\frac{XYZ}{1-XYZ}$&$p^{2r}T^{2r+1}$&$T^{-2}$&$1$ \\
1&2&$\frac{X^2YZ^2}{(1-XYZ)(1-XZ)}$&$p^{2r+1}T^{2r+1}$&$p^{-1}T^{-2}$&$1$ \\
2&2&$\frac{X^2Y^2Z}{(1-XYZ)(1-XY)}$&$p^{2r+1}T^{2r-1}$&$1$&$p^{-1}$ \\
3&3&$\frac{X^2YZ(1-X^2YZ)}{(1-XYZ)(1-XY)(1-XZ)(1-X)}$&$p^{2r+2}T^{2r+1}$&$p^{-1}T^{-2}$&$p^{-1}$ 
\end{tabular}
\caption{Computing $\Wsmptoff$ - lattices of type $(p^s,1,1)$.}
\label{smpt-table}
\end{table}

\noindent We get an expression for the desired generating function counting over lattices of type~$(p^s,1,1)$ in the boundary of this sector-family by substituting the respective Laurent monomials into the $F_j(X,Y,Z)$ and summing them up. More precisely:

\begin{eqnarray}
\lefteqn{\sum_{[\Lcen]}p^{w([\Lcen])\cdot2r}T^{w'([\Lcen])}=}\nonumber\\
&&\sum_{j=0}^3\sum_{(a,b,c)\in N_j}|\phi^{-1}(a,b,c)|\psi(a,b,c)= \nonumber \\
&&\sum_{j=0}^3(1-p^{-1})^{n_j-1}F_j(X,Y,Z)|_{\substack{  X=m_{jX}(p,T)\\  Y=m_{jY}(p,T)\\  Z=m_{jZ}(p,T)}}=\nonumber\\
&&\frac{p^{2r}T^{2r-1}(1-p^{2r+1}T^{2r+1})}{(1-p^{2r+1}T^{2r-1})(1-p^{2r+2}T^{2r+1})}, \label{gen.fu.smpt}
\end{eqnarray}
where~$[\Lcen]$ in the first sum ranges over all classes of lattices whose maximal element is of type~$(p^s,1,1),\,s\geq 1$ in the boundary of~$\SF$.

The generating function counting over lattices of type $(p^t,p^t,1)$, $t\geq1$, is clearly the same as in the previous case, i.e. is given by
\begin{equation}
\sum_{t\geq 1}p^{2t-2}\cdot p^{4rt}T^{(2r+2)t}=\frac{p^{4r}T^{2r+2}}{1-p^{4r+2}T^{2r+2}}.\label{gen.fu.offline}
\end{equation}

Counting over the interior of the fixed sector-family, i.e. over lattices of type~$(p^{s+t},p^t,1)$, $s,t\geq1$, is now easy. The generating function equals~$p$ times the product of the respective generating functions counting over the boundary. Indeed, there are~$p$ lattices of type~$(p^2,p^1,1)$ in~$\SF$, all of which carry the same weight~$w'([\Lcen])=4r+1$ and
\begin{equation*}
w'=w'|_{s=0}+w'_{t=0} \label{independence condition}
\end{equation*}
on the interior of~$\SF$. Therefore the parameters ~$s,t$ ``grow independently'' and the generating function is given by the product of the ones dealing with the cases~$t=0$ and~$s=0$, respectively. Thus, by~(\ref{gen.fu.smpt}) and~(\ref{gen.fu.offline}) we have
\begin{eqnarray}
\lefteqn{\sum_{[\Lcen]}p^{w([\Lcen])\cdot2r}T^{w'([\Lcen])}=}\nonumber\\
&&p\cdot\frac{p^{2r}T^{2r-1}(1-p^{2r+1}T^{2r+1})}{(1-p^{2r+1}T^{2r-1})(1-p^{2r+2}T^{2r+1})}\cdot\frac{p^{4r}T^{2r+2}}{1-p^{4r+2}T^{2r+2}},\label{gen.fu.smpt-offline interior}
\end{eqnarray}
where~$[\Lcen]$ in the first sum ranges over all classes of lattices whose maximal element is of type~$(p^{s+t},p^t,1),\,s,t\geq 1$ in the interior of~$\SF$.
Combining~(\ref{gen.fu.smpt}), (\ref{gen.fu.offline}) and (\ref{gen.fu.smpt-offline interior}) we get 
\begin{eqnarray}
\lefteqn{\Wsmptoff=}\nonumber\\
&&\frac{1}{\pthreetwo \ptwoone}+\frac{1}{\ptwoone}\left(\frac{p^{2r}T^{2r-1}(1-p^{2r+1}T^{2r+1})}{(1-p^{2r+1}T^{2r-1})(1-p^{2r+2}T^{2r+1})}+\frac{p^{4r}T^{2r+2}}{1-p^{4r+2}T^{2r+2}}\right)\nonumber \\
&&+\frac{p^{6r+1}T^{4r+1}(1-p^{2r+1}T^{2r+1})}{(1-p^{4r+2}T^{2r+2})(1-p^{2r+2}T^{2r+1})(1-p^{2r+1}T^{2r-1})}.\label{A_smpt/off}
\end{eqnarray} 

By simply counting the respective occurrences yields an expression for the generating $A(p,T)$:
\begin{eqnarray*}
A(p,T)&=&\left(\pthreetwo-|C(\Fp)|\right)\ptwoone \Woffoff+\\
&&\quad|C(\Fp)|\ptwoone \Wsmptoff\\
&=&A_1(p,T)+|C(\Fp)|A_2(p,T),
\end{eqnarray*}
where
\begin{eqnarray*}
A_1(p,T)&=&\frac{1+p^{2r}T^{2r+1}+p^{2r+1}T^{2r+1}+p^{4r}T^{2r+2}+p^{4r+1}T^{2r+2}+p^{6r+1}T^{4r+3}}{(1-p^{4r+2}T^{2r+2})(1-p^{2r+2}T^{2r+1})},\nonumber\\
A_2(p,T)&=&\frac{(1-T)(1+T)p^{2r}T^{2r-1}(1+p^{4r+1}T^{2r+2})}{(1-p^{4r+2}T^{2r+2})(1-p^{2r+2}T^{2r+1})(1-p^{2r+1}T^{2r-1})}.\nonumber\\
\end{eqnarray*}
But recall that in Lemma~\ref{lemma 1} in Section~\ref{section2} we described the local normal zeta function as a product of $A(p,p^{-s})$ and Riemann zeta functions. Thus we have produced the promised explicit formulae for the rational functions~$W_i(X,Y)$, $i=1,2$, and proved Theorem~\ref{main theorem}.

\verylongpage
\begin{proofnodot} (of Corollary~\ref{corollary to main theorem: fun. eq.}). Let $V$ be a non-singular, absolutely irreducible projective variety over $\Fp$ of dimension $n$. If $b_{V,e}$, $e\geq1$, denotes the number of $\mathbb{F}_{p^e}$-rational points of $V$, it is a well-known consequence of the {\sl rationality} of the Weil zeta function
$$Z_V(u)=\exp\left(\sum_{e=1}^\infty \frac{b_{V,e}u^e}{e}\right)$$
that there are complex numbers $\beta_{r,j}$, $r=0,\dots,2n$, $j=1,\dots,B_r$, $B_r\in \N$, such that 
$$b_{V,e}=\sum_{r=0}^{2n}(-1)^r\sum_{j=1}^{t_r}\beta^e_{r,j},$$
and that the function
\begin{eqnarray*}
\N_{\geq0}&\rightarrow&\N \\
e&\mapsto&b_{V,e}
\end{eqnarray*}
 has a {\sl unique} extension to $\Z$ (cf~\cite{DenefMeuser/91}, Lemma~2). 
The {\sl functional equation} of the Weil zeta function
$$Z_V(1/p^nu)=\pm (p^{n/2}u)^\chi Z_V(u),$$
where $\chi = \sum_{i=1}^{2n} (-1)^iB_i$, implies the $1-1$-correspondences
$$\left\{\frac{p^n}{\beta_{r,j}}|\;1\leq j \leq
B_r\right\}\stackrel{1-1}{\longleftrightarrow}\left\{\beta_{2n-r,i}|\;1
\leq i \leq B_{2n-r}\right\}$$
for $0\leq j \leq 2n$ (cf~\cite{Igusa/00}, p.~213). This gives 
$$b_{V,-e}=p^{-en}\,b_{V,e}$$ formally.
The corollary follows immediately from Observation~\ref{observation:fu.eq.} if we set $V=C$, $e=1$ and observe that $n=1$, $|C(\Fp)|=b_{C,1}$ and 
\begin{eqnarray*}
A_1(p,T)|_{\substack{p \rightarrow p^{-1} \\ T \rightarrow T^{-1}}}=p^3A_1(p,T)\\
A_2(p,T)|_{\substack{p \rightarrow p^{-1} \\
                    T \rightarrow T^{-1}}}=p^4A_2(p,T).
\end{eqnarray*}
\end{proofnodot}

\begin{example} In~\cite{duS-ecII/01}, du Sautoy gave an example of a \T-group $G(E)$ which is not finitely uniform. The associated Lie ring $\LieL(G(E))$ was presented as in Theorem~\ref{main theorem} with
$$R({\bf y})=\left(\begin{array}{ccc}
                            Dy_3 &y_1&y_2\\
                          y_1&y_3&0 \\
                          y_2&0&y_1 
                        \end{array}\right).$$

\noindent Note that $det(R({\bf y}))=Dy_1y_3^2-y_1^3-y_2y_3$, a polynomial defining the (projective) elliptic curve
$$E=y^2+x^3-Dx.$$
It follows readily from our Theorem~\ref{main theorem} that for almost all primes~$p$ we have
\begin{equation*}
\z_{G,p}(s)=W_1(p,p^{-s})+|E(\Fp)|W_2(p,p^{-s}),
\end{equation*}
where
\begin{eqnarray*}
W_1(X,Y)&=&\frac{(1+X^{6}Y^{7}+X^{7}Y^{7}+X^{12}Y^{8}+X^{13}Y^{8}+X^{19}Y^{15})}{\prod_{i=0}^6(1-X^iY)\cdot(1-X^9Y^{18})(1-^{14}Y^8)(1-X^8Y^7)},  \label{W1-ec} \\
W_2(X,Y)&=&\frac{(1-Y)(1+Y)X^6Y^5(1+X^{13}Y^8)}{\prod_{i=0}^6(1-X^iY)\cdot(1-X^9Y^{18})(1-X^{14}Y^8)(1-X^8Y^7)(1-X^7Y^5)}.  \label{W2-ec}
\end{eqnarray*}
The functional equation~(\ref{fun. eq. main theorem}) follows from the well-known fact that for each fixed curve $E$ there are complex numbers $\pi_p$ such that $\pi_p\cdot\overline{\pi_p}=p$ and 
$$|E(\mathbb{F}_{p^e})|=1-\pi_p^e-\overline{\pi_p}^e+p^e$$
and thus 
$$|E(\mathbb{F}_{p^{-e}})|=p^{-e}|E(\mathbb{F}_{p^e})|.$$

\end{example}

\bibliographystyle{amsplain} 
\bibliography{thebibliography}

\end{document}